\documentclass[11pt,reqno]{amsart}

\usepackage{lmodern}
\usepackage[T1]{fontenc}
\usepackage[utf8]{inputenc}
\usepackage[english]{babel}
\usepackage{microtype}
\usepackage{amsmath,amssymb,amsthm,mathtools,mathdots}
\usepackage{booktabs,enumerate}
\usepackage{xcolor}
\usepackage{hyperref}
\usepackage{algorithm2e}
\newtheorem{theorem}{Theorem}[section]

\theoremstyle{definition}

\newtheorem{example}[theorem]{Example}
\newtheorem{remark}[theorem]{Remark}

\newcommand{\setZ}{\mathbb{Z}}
\newcommand{\setQ}{\mathbb{Q}}
\newcommand{\setR}{\mathbb{R}}

\newcommand{\GT}{\mathrm{GT}}

\newcommand{\lam}{\lambda}
\newcommand{\wt}{\mathbf{w}}
\newcommand{\Ehr}{L}
\newcommand{\Ehrstar}{L^{*}}
\newcommand{\hvec}{\mathbf{h}^{*}}
\DeclareMathOperator{\interior}{int}

\newcommand{\defin}[1]{%
\relax\ifmmode%
\textcolor{blue}{#1}%
\else \textcolor{blue}{\emph{#1}}%
\fi%
}

\title[Ehrhart polynomials of GT polytopes via reciprocity]{Fast computation of Ehrhart polynomials of Gelfand--Tsetlin polytopes\\via Macdonald reciprocity}
\author{Per Alexandersson}
\address{Department of Mathematics, Stockholm University}
\email{per.w.alexandersson@gmail.com}
\date{}

\begin{document}

\begin{abstract}
We describe an efficient method for computing the Ehrhart polynomial of
Gelfand--Tsetlin polytopes arising from Kostka coefficients.
The key idea is to exploit Ehrhart--Macdonald reciprocity: evaluating
the Ehrhart polynomial at negative integers reduces to counting
\emph{strict} Gelfand--Tsetlin patterns, which are often zero or very
small for low dilations.  Combined with an adaptive strategy that chooses
the cheapest evaluation point (positive or negative) at each step, this
yields substantial practical speedups compared to general-purpose
polytope software.  We benchmark against $\mathtt{OSCAR}$/$\mathtt{polymake}$,
and illustrate the broader applicability of the method through
order polytopes and permutation posets.  The implementation is available in
the Rust \texttt{kostka} package, with related optimizations also incorporated
in the new \texttt{lrcalc-rs} replacement for \texttt{lrcalc}.
\end{abstract}

\maketitle

\section{Introduction}\label{sec:introduction}

Let $\lam/\mu$ be a skew shape and $\wt = (w_1, \dotsc, w_k)$ a weight
vector with $|\wt| = |\lam/\mu|$.  The \defin{Gelfand--Tsetlin polytope}
$\GT(\lam/\mu, \wt)$ is the set of all real-valued GT patterns of shape
$\lam/\mu$ and weight $\wt$.  Its integer points are in bijection with
semistandard Young tableaux of shape $\lam/\mu$ and weight $\wt$, so
\begin{equation}\label{eq:kostka}
  K_{\lam/\mu, \wt} = \#(\GT(\lam/\mu, \wt) \cap \setZ^d),
\end{equation}
where $d = \dim \GT(\lam/\mu, \wt)$.

By a theorem of E.~Rassart~\cite{Rassart2004}, the function
$n \mapsto K_{n\lam/n\mu, n\wt}$ is a polynomial in~$n$: the
\defin{Ehrhart polynomial} $\Ehr_P(n)$ of the polytope $P = \GT(\lam/\mu, \wt)$.
Computing this polynomial is important for understanding the structure of
Kostka coefficients and the geometry of GT polytopes.

The polynomiality of stretched Kostka and Littlewood--Richardson
coefficients, and the structure of the resulting Ehrhart polynomials,
have been studied by R.~King, C.~Tollu, and
F.~Toumazet~\cite{KingTolluToumazet2004,KingTolluToumazet2006LR}.
Their later hive-model paper~\cite{KingTolluToumazet2006LR}
computes several explicit stretched LR polynomials as illustrative
examples: the degree-8 polynomial for
$(\lam,\mu,\nu)=((4,3,3,2,1),(4,3,2,2,1),(7,4,4,4,3,2,1))$;
a factorization example with
$(\lam,\mu,\nu)=((9,7,6,2,0),(13,5,3,1,0),(14,12,11,5,4))$;
and degree-bound examples including
$(\lam,\mu,\nu)=((7,6,5,4,0,0),(7,7,7,4,0,0),(12,8,8,7,6,4,2))$
and
$(\lam,\mu,\nu)=((9,7,3,0,0),(9,9,3,2,0),(10,9,9,8,6))$.
These examples support conjectures on factorization, degree bounds,
and linear factors; they are not a systematic benchmark computation
of the family of GT/Kostka polytopes considered here.
Nzeutchap, Toumazet, and Butelle~\cite{NzeutchapToumazetButelle2006}
later developed a distributed method for stretched Kostka and
Littlewood--Richardson coefficients based on the hive model.  Their
approach computes values at positive dilation factors across a local
network and interpolates the resulting stretched polynomial; it also
uses hive inequalities to obtain an improved upper bound, denoted
\texttt{maxDeg}, on the degree before the computation begins.  Their
reported LR examples are arbitrary selected triples, with \texttt{maxDeg}
between 4 and 14; in several cases the real degree is strictly smaller,
including examples of real degree~0.
In contrast, for the GT/Kostka family considered here we compute the
dimension of the relevant GT polytope directly, so the exact degree of
the stretched Kostka polynomial is known before interpolation.  The
reciprocity step then reduces the number and cost of the needed
evaluations, since many negative evaluations are strict pattern counts
and often vanish.

Kostka coefficients arise as special cases of Littlewood--Richardson (LR) coefficients.
By the identity $K_{\lam/\nu,\mu} = \langle s_{\lam/\nu}, h_\mu \rangle
= \langle s_\lam, s_\nu h_\mu \rangle$, we construct a skew shape
$\rho/\sigma$ that is a disjoint union (no shared rows or columns) of
$\nu$ and single rows of lengths $\mu_1, \mu_2, \dotsc$, giving
$s_{\rho/\sigma} = s_\nu h_\mu$ and hence
$K_{\lam/\nu,\mu} = c^{\rho}_{\lam, \sigma}$.
Thus software for LR coefficients such as A.~Buch's
\texttt{lrcalc}~\cite{Buch2000} can be used for individual Kostka
evaluations.  However, the dilated partitions $n\lam$, $n\wt$
grow with~$n$, making LR computation expensive at large dilations.
Even \texttt{lrcalc} (highly optimized C code) takes over 8~minutes
for the 11~evaluations needed for $\GT((3,3,3), (1^9))$ at
dimension~10, compared to 1\,ms for our dynamic programming (DP)
approach with reciprocity
(Table~\ref{tab:lrcalc}).  Our horizontal-strip DP exploits the
specific structure of Kostka coefficients directly, without the
overhead of the LR reduction.

The standard approach evaluates $\Ehr_P(n)$ at $n = 1, 2, \dotsc, d+1$
using a DP algorithm for Kostka coefficients, then
recovers the polynomial by Lagrange interpolation.  Each evaluation
requires a DP that can be expensive for large dilations.

Our approach uses \defin{Ehrhart--Macdonald reciprocity} to access the
negative side of the polynomial for free or cheaply.  The same idea
also applies in other settings where one can count both lattice points
and relative interior lattice points efficiently; in particular, we
will later discuss order polytopes and an application to permutation
posets.

\section{Background}\label{sec:background}

\subsection{Ehrhart--Macdonald reciprocity}\label{subs:reciprocity}

For a rational polytope $P \subset \setR^d$ of dimension $d$ (one whose
vertices have rational coordinates), the function
$\Ehr_P(n) = \#(nP \cap \setZ^d)$ is a quasi-polynomial in~$n$, and a
polynomial when $P$ is a lattice polytope.  GT polytopes are rational but in general not lattice polytopes,
so one would \emph{a priori} expect $\Ehr_P$ to be only a
quasi-polynomial.  However, by Rassart's theorem~\cite{Rassart2004},
$\Ehr_P$ is a genuine polynomial---a case of so-called \emph{period collapse}.
Ehrhart--Macdonald reciprocity~\cite{Macdonald1971,BeckRobins2015}
gives:
\begin{equation}\label{eq:reciprocity}
  (-1)^d \Ehr_P(-n) = \defin{\Ehrstar_P(n)} \coloneqq \#(\interior(nP) \cap \setZ^d),
\end{equation}
where \defin{$\interior(nP)$} denotes the relative interior of $nP$.

\subsection{Gelfand--Tsetlin polytopes}\label{subs:GTpolytopes}

Recall that an SSYT of skew shape
$\lam/\mu$ and content $\wt = (w_1, \dotsc, w_k)$ corresponds
bijectively to a chain of partitions
\[
  \mu = \alpha^{(0)} \subseteq \alpha^{(1)} \subseteq \dotsb
  \subseteq \alpha^{(k)} = \lam,
\]
where each $\alpha^{(i)} / \alpha^{(i-1)}$ is a horizontal strip
of size~$w_i$.  The \defin{Gelfand--Tsetlin polytope}
$\GT(\lam/\mu, \wt) \subset \setR^m$ is the convex polytope obtained
by relaxing the integrality constraint: the entries of the
intermediate partitions $\alpha^{(1)}, \dotsc, \alpha^{(k-1)}$ are
allowed to be real, subject to the interlacing inequalities
\[
  \alpha^{(i-1)}_j \leq \alpha^{(i)}_j \leq \alpha^{(i-1)}_{j-1}
\]
and the row-sum constraints
$|\alpha^{(i)}| - |\alpha^{(i-1)}| = w_i$.
The lattice points of $\GT(\lam/\mu, \wt)$ are precisely the SSYT,
so $K_{\lam/\mu, \wt} = \#(\GT(\lam/\mu, \wt) \cap \setZ^m)$.

\subsection{Strict GT patterns}\label{subs:strictGT}

The interior lattice points of $n \cdot \GT(\lam/\mu, \wt)$ are
precisely the \defin{strict GT patterns}: those where every interlacing
inequality that is \emph{not} forced as an equality in the affine hull
of the polytope is strict.  Equivalently, one imposes strictness only on
the non-globally-tight inequalities.  Thus
\begin{equation}\label{eq:strict-kostka}
  \Ehrstar_P(n) = K^{\mathrm{strict}}_{n\lam/n\mu, n\wt},
\end{equation}
where $K^{\mathrm{strict}}$ counts strict semistandard tableaux
(equivalently, GT patterns with strictness on all non-forced internal
inequalities).  This is computed by the same DP as the ordinary Kostka
coefficient, but with strict inequality constraints.

\section{The adaptive reciprocity algorithm}\label{sec:algorithm}

To compute $\Ehr_P(n)$ for a polytope of dimension~$d$, we need $d+1$
evaluations at distinct integers.  We always have $\Ehr_P(0) = 1$
for free, and by~\eqref{eq:reciprocity} and~\eqref{eq:strict-kostka},
each negative evaluation $(-1)^d \Ehr_P(-n) = K^{\mathrm{strict}}_{n\lam, n\mu, n\wt}$
reduces to a strict Kostka computation.  The key observation is that
strict Kostka coefficients are often \emph{zero} for small~$n$:
if any part of $\wt$ is~1, then strict GT patterns require the
corresponding row to have distinct entries in a range that may be
empty for small dilations.  Even when the strict count is nonzero,
its DP state space is typically much smaller than the ordinary
Kostka DP for the same dilation.

This motivates an adaptive strategy:

\begin{algorithm}[H]
\SetAlgoLined
\KwIn{Shape $\lam/\mu$, weight $\wt$, dimension $d$}
\KwOut{Ehrhart polynomial $\Ehr_P(n)$}
Initialize: $\text{points} \gets \{(0, 1)\}$\;
$p \gets 0$, $q \gets 0$\;
$\text{cost}_+ \gets 1$, $\text{cost}_- \gets 1$\;
\While{$|\text{points}| \leq d$}{
  \eIf{$\text{cost}_- \leq \text{cost}_+$}{
    $q \gets q + 1$\;
    $v \gets K^{\mathrm{strict}}_{q \lam, q \mu, q \wt}$\;
    Append $(-q, \; (-1)^d v)$ to points\;
    $\text{cost}_- \gets v$\;
  }{
    $p \gets p + 1$\;
    $v \gets K_{p \lam, p \mu, p \wt}$\;
    Append $(p, v)$ to points\;
    $\text{cost}_+ \gets v$\;
  }
}
\Return Lagrange interpolation of points\;
\caption{Adaptive Ehrhart polynomial computation}
\end{algorithm}

The ``cost'' heuristic uses the most recent Kostka/strict-Kostka value
as a proxy for the DP state-space size at the next dilation.  The
negative side typically wins for the first several points because
strict Kostka equals zero.

\begin{remark}\label{rem:freePoints}
The negative evaluations are completely free when $K^{\mathrm{strict}} = 0$.
For a weight vector with $k$ ones, the strict Kostka coefficient
$K^{\mathrm{strict}}_{n\lam, n\mu, n\wt}$ is zero whenever $n$ is
small enough that the strict interlacing constraints for those
unit-weight rows cannot be satisfied.  In practice, half or more of
the $d+1$ required evaluation points are obtained this way.
\end{remark}

\section{Implementation}\label{sec:implementation}

Our implementation is in Rust.  The specialized research code used for
the computations in this article is available as the \texttt{kostka}
crate.\footnote{\url{https://github.com/PerAlexandersson/kostka}}
It contains the GT/Kostka dynamic programs, the Ehrhart interpolation and
reciprocity routines, and the strict-count variants used below.
The order-polytope experiments discussed later use the same
reciprocity/interpolation viewpoint, but are conceptually separate from
the \texttt{lrcalc} compatibility layer.
At its core is a DP algorithm for (skew) Kostka coefficients
$K_{\lam/\mu, \wt}$ based on horizontal-strip removal, using
\texttt{BigUint} arithmetic for exact computation.
A variant with strict inequality constraints in the DP transitions
handles the strict Kostka coefficients
$K^{\mathrm{strict}}_{\lam/\mu, \wt}$.
The dimension of $\GT(\lam/\mu, \wt)$ is determined by counting
free entries in the GT pattern with propagation of forced values.
This gives the exact degree of the stretched Kostka polynomial, rather
than only an a priori upper bound; see J.~A.~De~Loera and
T.~B.~McAllister~\cite{DeLoeraMcAllister2006} and
Nzeutchap, Toumazet, and Butelle~\cite{NzeutchapToumazetButelle2006}
for related degree computations.
The Ehrhart polynomial is then recovered from evaluation points
via Lagrange interpolation over~$\setQ$.

The \texttt{kostka} CLI supports:
\begin{verbatim}
  kostka ehrhart --lambda 3,2,1 -w 1,1,1,1,1,1
  kostka ehrhart --lambda 4,3 --mu 2,1 -w 2,2,2
  kostka hstar --lambda 4,3,2,1 -w 1,1,1,1,1,1,1,1,1,1
  kostka hstar --lambda 4,3,1 --mu 2,1 -w 2,2,1
\end{verbatim}
The tool also computes Kostka coefficients associated with flagged
Schur functions, Littlewood--Richardson coefficients, and batch
computations over all weights for a given shape.  We refer to the
Git repository for full documentation.

There is also a separate Rust implementation aimed at the classical
\texttt{lrcalc} interface:
\texttt{lrcalc-rs}.\footnote{%
\url{https://github.com/PerAlexandersson/lrcalc-rs}}
That project is not the specialized Ehrhart engine used in
Tables~\ref{tab:benchmarks} and~\ref{tab:lrcalc}; rather, it is a
drop-in-oriented implementation of the \texttt{liblrcalc} ABI and
command-line tools, with native code for Littlewood--Richardson
coefficients, Schur products, skew Schur expansions, Schubert products,
and fast Kostka special cases.  The speedups in the present article
come from exploiting the GT/Kostka and reciprocity structure directly,
while \texttt{lrcalc-rs} packages related optimizations in a broader
LR/Schur-function compatibility library.

\section{Benchmarks}\label{sec:benchmarks}

We compare our method against two alternatives:
the $\mathtt{OSCAR}$/$\mathtt{polymake}$ pipeline~\cite{GawrilowJoswig2000,OskarThesis}
and A.~Buch's \texttt{lrcalc}~\cite{Buch2000}, which reduces Kostka
evaluations to LR coefficient computations.
($\mathtt{OSCAR}$ is a Julia framework that calls $\mathtt{polymake}$ for polyhedral
computations, which in turn delegates lattice-point enumeration
and $h^*$-vector computation to $\mathtt{Normaliz}$~\cite{BrunsIchim2010};
thus ``$\mathtt{OSCAR}$,'' ``$\mathtt{polymake}$,'' and ``$\mathtt{Normaliz}$'' refer to different
layers of the same computational stack.)
Tables~\ref{tab:benchmarks} and~\ref{tab:lrcalc} summarize the results.

The $\mathtt{OSCAR}$/$\mathtt{polymake}$ approach requires constructing the full vertex
representation of the GT polytope---expensive in high dimensions---and
does not exploit Ehrhart--Macdonald reciprocity.  Our method avoids
vertex enumeration entirely: each evaluation of $\Ehr_P(n)$ or
$\Ehrstar_P(n)$ reduces to a Kostka coefficient DP whose complexity
depends on the dilation~$n$ (not the number of vertices), and
reciprocity provides many free evaluation points where strict Kostka
equals zero.

\begin{table}[!ht]
\centering
\small
\caption{%
  Comparison of Ehrhart polynomial computation times.
  ``$d$'' is the polytope dimension.
Column ``\texttt{kostka}'' is our DP + reciprocity method (desktop).
Column ``$\mathtt{OSCAR}$'' is from~\cite{OskarThesis} using $\mathtt{OSCAR}$/$\mathtt{polymake}$ (laptop;
  because the hardware differs, the reported ratios should be read as
  indicative rather than exact speedups).
  Cases marked ``---'' exceeded the available time/memory.
}
\label{tab:benchmarks}
\begin{tabular*}{\textwidth}{@{\extracolsep{\fill}}llrrrr@{}}
\toprule
$\lam$ & $\wt$ & $d$ & \texttt{kostka} & $\mathtt{OSCAR}$ & speedup \\
\midrule
\multicolumn{6}{@{}l}{\emph{Direct comparisons (partitions of 10, same cases):}} \\[2pt]
$(8,1,1)$       & $(2,1^8)$  & 12 & 2\,ms   & 666\,s   & $3 \times 10^5$ \\
$(4,4,2)$       & $(2,1^8)$  & 12 & 2\,ms   & 182\,s   & $9 \times 10^4$ \\
$(4,3,3)$       & $(2,1^8)$  & 12 & 2\,ms   & 230\,s   & $10^5$  \\
$(5,3,2)$       & $(2,1^8)$  & 13 & 4\,ms   & 764\,s   & $2 \times 10^5$ \\
$(7,2,1)$       & $(2,1^8)$  & 13 & 3\,ms   & 1347\,s  & $4 \times 10^5$ \\
$(7,1,1,1)$     & $(2,1^8)$  & 15 & 5\,ms   & 467\,s   & $9 \times 10^4$ \\
$(4,2,2,2)$     & $(2,1^8)$  & 15 & 9\,ms   & 125\,s   & $10^4$  \\
$(3,3,3,1)$     & $(2,1^8)$  & 15 & 8\,ms   & 573\,s   & $7 \times 10^4$ \\
$(5,2,2,1)$     & $(2,1^8)$  & 17 & 21\,ms  & 1988\,s  & $9 \times 10^4$ \\
$(5,3,1,1)$     & $(2,1^8)$  & 17 & 18\,ms  & 1179\,s  & $7 \times 10^4$ \\
\midrule
\multicolumn{6}{@{}l}{\emph{Cases beyond $\mathtt{OSCAR}$'s reach:}} \\[2pt]
$(4,3,2,1)$     & $(1^{10})$  & 21 & 104\,ms  & ---  & \\
$(5,5,5)$       & $(1^{15})$  & 22 & 45\,ms   & ---  & \\
$(5,3,3,1,1,1)$ & $(2^4,1^6)$ & 26 & 17\,s    & ---  & \\
\bottomrule
\end{tabular*}
\end{table}

\begin{remark}\label{rem:speedup}
The observed wall-clock ratios are typically between $10^4$ and $10^5$.
The $\mathtt{OSCAR}$/$\mathtt{polymake}$ timings are from~\cite{OskarThesis},
computed on a laptop using the Julia $\mathtt{OSCAR}$ framework
(which calls $\mathtt{polymake}$/$\mathtt{Normaliz}$ internally for $h^*$-vector computation).
The total computation time reported in~\cite{OskarThesis} for all
partitions of size~10 was over 412~hours.  Our method computes the same
cases in under 1~second total, and extends to much larger cases
(55{,}000 polytopes up to dimension~26 in our database).

For standard weight $\wt = (1^N)$, the reciprocity method is especially
effective: $K^{\mathrm{strict}}_{n\lam, n\wt} = 0$ for all small $n$,
so many evaluation points are free.  For example, computing the
degree-21 Ehrhart polynomial for $\GT((4,3,2,1), (1^{10}))$ requires
22 evaluation points, of which the first 9 negative points all give
$K^{\mathrm{strict}} = 0$.
\end{remark}

\begin{table}[!ht]
\centering
\caption{%
  Comparison with \texttt{lrcalc}~\cite{Buch2000}.
  Each Kostka evaluation $K_{n\lam, n\wt}$ is computed via the
  Stanley reduction to an LR coefficient $c^{n\rho}_{n\lam, n\sigma}$
  and then calling \texttt{lrcalc lrcoef}.  Both methods run on the
  same desktop.  The \texttt{lrcalc} column is the total time for
  all $d+1$ positive-dilation evaluations (no reciprocity).
}
\label{tab:lrcalc}
\begin{tabular*}{\textwidth}{@{\extracolsep{\fill}}llrrrr@{}}
\toprule
$\lam$ & $\wt$ & $d$ & \texttt{kostka} & \texttt{lrcalc} & speedup \\
\midrule
$(3,2,1)$   & $(1^6)$  &  7 & ${<}\,1$\,ms & 86\,ms    & $\sim\!10^2$ \\
$(3,3,3)$   & $(1^9)$  & 10 & 1\,ms         & 491\,s    & $\sim\!6 \times 10^3$ \\
$(4,3,2)$   & $(1^9)$  & 13 & 4\,ms         & $>\!$10\,min & $>\!10^5$ \\
\bottomrule
\end{tabular*}
\end{table}

\begin{remark}\label{rem:lrcalcSlowdown}
The slowdown of \texttt{lrcalc} at higher dilations is expected:
computing $c^{n\rho}_{n\lam, n\sigma}$ requires enumerating
Littlewood--Richardson tableaux of shape $n\rho / n\lam$ with
content $n\sigma$, whose count grows rapidly with~$n$.  Our DP, by
contrast, processes one row of the GT pattern at a time and benefits
from the horizontal-strip structure.  The reciprocity advantage
(free evaluations at negative integers) further widens the gap.
The distributed hive computation of
Nzeutchap--Toumazet--Butelle~\cite{NzeutchapToumazetButelle2006} is
closer in spirit to interpolation: it parallelizes the computation of
positive dilation values for stretched LR coefficients.  Our comparison
here is different in two respects: we specialize to the GT/Kostka
setting, where the horizontal-strip DP is very small, and we use
Ehrhart--Macdonald reciprocity to replace many positive evaluations by
strict negative-side evaluations.
\end{remark}

\section{The \texorpdfstring{$h^{*}$}{h*}-vector}\label{sec:hstar}

For a \emph{rational} polytope $P$ with denominator $q$
(the smallest positive integer such that $qP$ is a lattice polytope),
the lattice-point counting function $\Ehr_P(n) = |nP \cap \setZ^d|$ is a
\defin{quasi-polynomial} of degree $d$ and period dividing $q$.
The Ehrhart series then satisfies
\begin{equation}\label{eq:hstar-rational}
  \sum_{n \geq 0} \Ehr_P(n) \, z^n
  = \frac{h^*_q(P;z)}{(1-z^q)^{d+1}},
\end{equation}
where $h^*_q(P;z)$ is a polynomial of degree at most $q(d+1)-1$
with non-negative integer coefficients,
see~\cite{BeckSottile2007} and~\cite{Stanley1993}.
Since GT polytopes have period collapse by Rassart's theorem,
we define the \defin{$h^{*}$-vector} $(h^*_0, \dotsc, h^*_d)$
as in the lattice polytope case ($q=1$):
\begin{equation}\label{eq:hstar-lattice}
  \sum_{n \geq 0} \Ehr_P(n) \, z^n
  = \frac{h^*_0 + h^*_1 z + \dotsb + h^*_d z^d}{(1-z)^{d+1}}.
\end{equation}
Comparing~\eqref{eq:hstar-rational} and~\eqref{eq:hstar-lattice}
via the factorization
$(1-z^q)^{d+1} = (1-z)^{d+1}(1+z+\dotsb+z^{q-1})^{d+1}$,
we see that
$h^*_q(P;z) = (h^*_0 + \dotsb + h^*_d z^d) \cdot (1+z+\dotsb+z^{q-1})^{d+1}$.
In particular, non-negativity of the $h^*_i$
implies non-negativity of the coefficients of $h^*_q(P;z)$.
R.~King, C.~Tollu, and F.~Toumazet~\cite{KingTolluToumazet2004}
conjectured that the Ehrhart polynomial $\Ehr_P(n)$ of a GT polytope
has non-negative coefficients (in the standard monomial basis).
Note that this is an independent property from $h^*$-positivity:
neither implies the other in general.
Our computations verify both properties for all straight-shape GT
polytopes $\GT(\lam,(1^{|\lam|}))$ with $|\lam| \leq 15$---over
54{,}000 cases---providing further evidence for the
conjecture. 
See also \cite{AlexanderssonAlhajjar2018,AlexanderssonOguz2023x}
for related positivity conjectures for Ehrhart polynomials
of polytopes related to key polynomials,
flagged Schur functions, and cylindric Schur functions.
For recent progress on Ehrhart positivity for marked order polytopes,
relevant to skew GT polytopes without a weight condition,
see~\cite{JochemkoMenon2026x}.

As observed by C.~Haase and T.~B.~McAllister~\cite{HaaseMcAllister2007},
there exist rational polytopes with period collapse whose Ehrhart
polynomial is not the Ehrhart polynomial of any lattice polytope.
It is natural to ask whether the same phenomenon occurs for GT
polytopes: is $\Ehr_{\GT(\lam/\mu, \wt)}(n)$ always realizable as
the Ehrhart polynomial of a lattice polytope, or are GT polytopes
genuinely different in this respect?

\section{Order polytopes}\label{sec:orderPolytopes}

The reciprocity approach extends naturally to \defin{order polytopes}.
For a finite poset~$P$ on $\{0,1,\dotsc,n-1\}$, the
\defin{order polytope}~\cite{Stanley1986}
\[
  \mathcal{O}(P) = \bigl\{x \in [0,1]^n : x_a \leq x_b \text{ for all } a <_P b \bigr\}
\]
has the property that its lattice points at dilation~$t$ are weakly
order-preserving maps $P \to \{0,1,\dotsc,t\}$, so
$\Ehr_{\mathcal{O}(P)}(t) = \Omega(P, t+1)$, where $\Omega(P,k)$
is Stanley's order polynomial.  Interior lattice points at dilation~$t$
correspond to strictly order-preserving maps $P \to \{1,\dotsc,t-1\}$,
giving Ehrhart--Macdonald reciprocity in the form
\begin{equation}\label{eq:order-reciprocity}
  \Ehr_{\mathcal{O}(P)}(-k)
  = (-1)^n \, \bar{\Omega}(P, k-1),
\end{equation}
where $\bar{\Omega}(P,m)$ counts strict order-preserving maps
$P \to \{1,\dotsc,m\}$.

\subsection{Frontier DP for order-preserving maps}\label{subs:frontierDP}

The standard approach to evaluating $\Omega(P,k)$ enumerates all
order-preserving maps by backtracking over vertices, giving
$O(k^n)$ worst-case time.  We introduce a \defin{frontier DP} that
processes vertices in a natural labeling order, maintaining as state
only the values of \emph{live} vertices---those with at least one
unprocessed child.  When a vertex~$v$ has no children among the
remaining vertices, its value does not affect future constraints,
so all valid choices are collapsed into a single multiplier
$(k - \ell_v + 1)$, where $\ell_v$ is the lower bound from parent
constraints.  States that agree on all live vertex values are merged.

The complexity is $O(n \cdot k^{w})$, where $w$ is the maximum
frontier width (the largest number of simultaneously live vertices).
For chains and zigzag (fence) posets, $w = 1$; for antichains,
$w = 0$ (each vertex is independent); for Young diagram posets of
shape $(\lambda_1, \dotsc, \lambda_r)$, $w \leq \lambda_1$.

\subsection{Order polytope benchmarks}\label{subs:orderBenchmarks}

We benchmark three methods for computing the $h^*$-vector of order polytopes:
\begin{enumerate}
\item \textbf{Rust DP}: our frontier DP with Ehrhart--Macdonald reciprocity
  and \texttt{BigRational} interpolation;
\item \textbf{linext}: Stanley's formula---enumerate all linear extensions
  $\mathcal{L}(P)$ and tabulate descents (streaming, $O(n)$ memory);
\item \textbf{$\mathtt{polymake}$}: the $\mathtt{polymake}$/$\mathtt{Normaliz}$ backend for general
  lattice-point enumeration.
\end{enumerate}

\begin{table}[!ht]
\centering
\small
\caption{%
  Comparison of order polytope $h^*$-vector computation times.
  All three methods run on the same machine (desktop, single thread).
  ``---'' indicates the method was not run
  (the linear extension method is infeasible
  for large fence posets due to the exponential number of extensions).
}
\label{tab:orderPolyBench}
\begin{tabular*}{\textwidth}{@{\extracolsep{\fill}}lrrrrr@{}}
\toprule
Poset & $n$ & Rust DP & linext & $\mathtt{polymake}$ \\
\midrule
Shape $(3,2,1)$    &  6 & 0.4\,ms & 2\,\textmu s  & 60\,ms  \\
Fence(8)           &  8 & 0.9\,ms & 50\,\textmu s & 90\,ms  \\
Shape $(4,3,2,1)$  & 10 & 2\,ms   & 30\,\textmu s & 90\,ms  \\
Shape $(3,3,3)$    &  9 & 1.3\,ms & 5\,\textmu s  & 60\,ms  \\
Fence(10)          & 10 & 1.8\,ms & 1.7\,ms       & 520\,ms \\
2-alt(10)          & 10 & 1.8\,ms & 1.7\,ms       & 490\,ms \\
\midrule
Shape $(5,4,3,2,1)$ & 15 & 9\,ms  & 11\,ms        & 2.5\,s  \\
Shape $(4,4,4,4)$  & 16 & 20\,ms  & 2\,ms         & 270\,ms \\
Shape $(5,5,5)$    & 15 & 20\,ms  & 0.4\,ms       & 120\,ms \\
Fence(14)          & 14 & 6\,ms   & ---            & $>\!5$\,s \\
\midrule
Shape $(5,5,5,5)$  & 20 & 175\,ms & ---            & $>\!5$\,s \\
Shape $(6,5,4,3,2,1)$ & 21 & 100\,ms & ---        & ---     \\
Fence(20)          & 20 & 23\,ms  & ---            & ---     \\
Shape $(6,6,6)$    & 18 & 230\,ms & ---            & ---     \\
Shape $(6,6,6,6)$  & 24 & 6.7\,s  & ---            & ---     \\
\bottomrule
\end{tabular*}
\end{table}

\begin{remark}\label{rem:orderPolySpeedup}
Table~\ref{tab:orderPolyBench} reveals that the three methods have
complementary strengths.  For \emph{few linear extensions}
(shapes with small hook lengths), Stanley's formula is fastest: the
staircase $(4,3,2,1)$ has only 768 linear extensions, so direct
enumeration takes 30\,\textmu s.
For \emph{narrow} posets (chains, fences), the frontier DP dominates
because $w=1$ gives polynomial time $O(nk)$, while the number of
linear extensions grows exponentially---fence(14) already exhausts
the 5\,s cap for $\mathtt{polymake}$, but the DP finishes in 6\,ms.
$\mathtt{polymake}$/$\mathtt{Normaliz}$ is a general-purpose tool that does not exploit
the order polytope structure; its triangulation-based approach
scales with the number of vertices of the polytope, not the
poset width.  For the largest cases ($n \geq 20$), only the
frontier DP is feasible.
\end{remark}

By Stanley's theorem~\cite[Thm.~4.5.14]{StanleyEC1},
for a naturally labeled poset~$P$, the $h^*$-vector of
$\mathcal{O}(P)$ equals the $P$-Eulerian polynomial:
\[
  h^*_i = \#\{\sigma \in \mathcal{L}(P) : \mathrm{des}(\sigma) = i\},
\]
where $\mathcal{L}(P)$ is the set of linear extensions and
$\mathrm{des}(\sigma)$ counts descents.
Counting linear extensions is $\#\mathsf{P}$-complete in
general~\cite{BrightwellWinkler1991}, and even for structured posets the
number $|\mathcal{L}(P)|$ can be enormous---e.g., $|\mathcal{L}(P)| = n!$
for the antichain on $n$~elements.  Our Ehrhart/reciprocity approach
avoids linear extension enumeration entirely: it requires only $\sim n/2$
evaluations of the order polynomial via the frontier DP, each costing
$O(n \cdot k^w)$.

\begin{example}\label{ex:fenceHstar}
For the fence (zigzag) poset on 10 elements:
\[
  \hvec = (1, 133, 2475, 12331, 20641, 12331, 2475, 133, 1).
\]
This is palindromic, as expected for a Gorenstein polytope.
\end{example}

\subsection{Non-real-rooted \texorpdfstring{$h^*$}{h*}-vectors from permutation posets}\label{subs:nonRR}

As an application of the fast $h^*$-computation, we investigate
real-rootedness of $h^*$-polynomials for \defin{permutation posets}.
Given $w \in S_n$, the permutation poset $P_w$ is the partial order
on $\{1,\dotsc,n\}$ with $i <_{P_w} j$ iff $i < j$ and $w_i < w_j$.
Since $P_w$ is always naturally labeled, $h^*(\mathcal{O}(P_w); t)$
equals the $P$-Eulerian polynomial by Stanley's theorem.

Starting from a $321$-avoiding permutation of Stembridge,
\[
  w_0 = (2,4,6,8,10,1,12,3,15,5,17,7,9,11,13,14,16) \in S_{17},
\]
whose permutation poset has
$h^*(t) = 1 + 32t + 336t^2 + 1420t^3 + 2534t^4 + 1946t^5 + 658t^6 + 86t^7 + 3t^8$
(ultra-log-concave but \emph{not} real-rooted),
we searched over 34{,}226 nearby $4321$-avoiding permutations
(within 3~arbitrary transpositions) and found exactly 22 with
non-real-rooted~$h^*$.  All 22 counterexamples are
ultra-log-concave, hence also log-concave.

\begin{example}\label{ex:321containing}
Almost all non-real-rooted examples are $321$-avoiding.
A notable exception, found at transposition distance~3 from~$w_0$, is
\[
  w = (3,4,6,8,10,\underline{12,2,1},15,5,17,7,9,11,13,14,16) \in S_{17},
\]
which \emph{contains} the pattern~$321$. 
 Its $h^*$-polynomial is
\[
  h^*(t) = 1 + 41t + 525t^2 + 2596t^3 + 5349t^4 + 4731t^5 + 1849t^6 + 284t^7 + 12t^8,
\]
with non-real roots at approximately $-2.15 \pm 0.041\,i$.
This is the only $321$-containing example among the 22 failures
at $n = 17$, and has the \emph{smallest} imaginary part of all
examples found---it is the closest to being real-rooted.
\end{example}

\begin{example}\label{ex:S28nonRR}
A larger example occurs for the permutation
\[
\begin{aligned}
  w ={}&(9,10,1,2,3,4,5,12,15,16,17,18,19,6,\\
        &7,8,11,20,21,22,23,13,25,26,27,28,14,24)
  \in S_{28}.
\end{aligned}
\]
The permutation poset~$P_w$ has
\[
\begin{aligned}
  h^*(t)={}&1 + 66t + 1500t^2 + 15582t^3 + 81644t^4
  + 223486t^5  \\
  &+ 320052t^6 + 232424t^7 + 77660t^8 + 8560t^9 .
\end{aligned}
\]
This polynomial is ultra-log-concave, but it is not real-rooted.
\end{example}

\section{Further directions}\label{sec:furtherDirections}

Our method applies whenever one has a combinatorial interpretation
of both $\Ehr_P(n)$ and the interior count $\Ehrstar_P(n)$ that
admits a DP or transfer-matrix algorithm, together with a regime
where $\Ehrstar_P(n) = 0$ for small~$n$.
GT polytopes are closely related to \defin{flow polytopes} on directed acyclic
graphs~\cite{LiuMeszarosDizier2019flow}.

Another natural test case is the Birkhoff polytope~$B_\ell$.
Here $\Ehr_{B_\ell}(n)$ and $\Ehrstar_{B_\ell}(n)$ count
nonnegative and positive integer magic squares with line sum~$n$,
so the reciprocity setup again applies.
For Birkhoff polytopes, however, the symmetry-based interpolation trick
is classical rather than new.  Beck and Pixton~\cite{BeckPixton2003}
record the standard consequence of Ehrhart--Macdonald reciprocity:
from reciprocity one gets
$\Ehr_{B_\ell}(-\ell-t)=(-1)^{\ell-1}\Ehr_{B_\ell}(t)$ and
$\Ehr_{B_\ell}(-1)=\cdots=\Ehr_{B_\ell}(-\ell+1)=0$.
They explicitly describe the resulting strategy: compute the first
$\binom{\ell-1}{2}$ values, use the symmetry and trivial zeros, and
interpolate the polynomial.  They also note that the volume of~$B_8$
had been computed by essentially this method, combined with additional
computational tricks.  Beck--Pixton's own main contribution in that
paper is a different residue-based method for computing $H_\ell(t)$
and the volumes of larger Birkhoff polytopes.
For context, our implementation computes the Ehrhart polynomial in under 10\,ms for
$B_3$ and~$B_4$, in 0.28\,s for~$B_5$, and in about 265\,s for~$B_6$.
Moreover, $\Ehrstar_{B_\ell}(n)=0$ for $1 \leq n < \ell$,
so the first $\ell-1$ negative evaluation points are free.
Thus the new contribution in this direction would be the GT-DP implementation
and its adaptation of that classical Birkhoff symmetry trick, not the trick
itself.
This suggests that the method is viable for smaller Birkhoff polytopes,
although the growth by $\ell=6$ is already steep.

Matroid polytopes are another natural family
where the reciprocity approach may be effective.

\subsection*{Acknowledgements}

Some Rust code, installation work, and benchmark scripts were developed with
assistance from OpenAI Codex.

\bibliographystyle{amsalpha}
\bibliography{ehrhart-reciprocity}

\end{document}